\documentclass{article} 
\usepackage{iclr2026_conference,times}
\usepackage{graphicx}

\usepackage[normalem]{ulem}
\UseRawInputEncoding

\usepackage{amsmath,amsfonts,bm}

\def\eqref#1{equation~\ref{#1}}

\def\1{\bm{1}}

\DeclareMathAlphabet{\mathsfit}{\encodingdefault}{\sfdefault}{m}{sl}
\SetMathAlphabet{\mathsfit}{bold}{\encodingdefault}{\sfdefault}{bx}{n}

\newcommand{\R}{\mathbb{R}}

\usepackage{hyperref}
\usepackage{url}
\iclrfinalcopy   

\title{Learning embeddings of non-linear PDEs: the Burgers' equation}
\author{Pedro Taranc\'on-\'Alvarez \\
Departament de F\'isica Qu\'antica i Astrof\'i­sica, Universitat de Barcelona \\
Institut de Ciencies del Cosmos (ICC), Universitat de Barcelona,\\
Marti i Franques 1,
ES-08028, \\Barcelona, Spain.\\
\texttt{\url{pedro.tarancon@fqa.ub.edu}}
\And
Leonid Sarieddine \\
Departament de F\'isica Qu\'antica i Astrof\'i­sica, Universitat de Barcelona \\
Institut de Ciencies del Cosmos (ICC), Universitat de Barcelona,\\
Marti i Franques 1,
ES-08028, \\Barcelona, Spain.\\
\texttt{\url{leonid.sarieddine@icc.ub.edu}}
\And
Pavlos Protopapas \\
Institute for Applied Computational Science \\
Harvard University \\
Cambridge, MA, USA \\
\texttt{\url{pprotopapas@g.harvard.edu}}
\AND
Raul Jimenez \\
Institut de Ciencies del Cosmos (ICC), Universitat de Barcelona,\\
Marti i Franques 1,
ES-08028, \\Barcelona, Spain.\\
Institucio Catalana de Recerca i Estudis Avancats (ICREA),\\
Passeig Lluis Companys 23, ES-08010,\\ 
Barcelona, Spain\\
\texttt{\url{raul.jimenez@icc.ub.edu}}
}

\begin{document}
\maketitle

\begin{abstract}
Embeddings provide low-dimensional representations that organize complex function spaces and support generalization. They provide a geometric representation that supports efficient retrieval, comparison, and generalization. In this work we generalize the concept to Physics Informed Neural Networks. We present a method to construct solution embedding spaces of nonlinear partial differential equations using a multi-head setup, and extract non-degenerate information from them using principal component analysis (PCA). We test this method by applying it to viscous Burgers' equation, which is solved simultaneously for a family of initial conditions and values of the viscosity. A shared network body learns a latent embedding of the solution space, while linear heads map this embedding to individual realizations. By enforcing orthogonality constraints on the heads, we obtain a principal-component decomposition of the latent space that is robust to training degeneracies and admits a direct physical interpretation. The obtained components for Burgers' equation exhibit rapid saturation, indicating that a small number of latent modes captures the dominant features of the dynamics.
\end{abstract}

\section{Introduction}
Scientific machine learning (ML) has long targeted partial differential equations (PDEs) regimes where solutions are nonlinear, stiff, and multiscale.
For these problems, it is valuable not only to predict accurately the solution, but also to learn \emph{low-dimensional coordinates} that organize families of solutions across initial/boundary data and physical parameters. Such coordinates could support fast surrogates, inverse inference, and model reduction.

This ``representation'' view is especially relevant since many operator-learning methods for solving differential equations (DEs) assume that the solution set concentrates near a low-dimensional manifold.
Here we make that assumption explicit and measurable and focus on the geometrical structure of this space.
We learn an embedding or latent space using a shared PINN body and then quantify its intrinsic dimensionality through a PCA of the latent-space covariance matrix.

As a concrete and well-controlled testbed, we consider the 1D viscous Burgers' equation. This is a nonlinear, second-order PDE that develops steep gradients and shock-like (viscosity-regularized) features as $\nu$ decreases. Burgers is also a standard surrogate for aspects of turbulent transport and provides a transparent testbed for methods aimed at more complex systems. The equation reads
\begin{equation}
\partial_t u + u\,\partial_x u = \nu\,\partial_{xx}u,
\qquad (x,t,\nu)\in[-5,5]\times[0,5] \times [10^{-2},1].
\label{eq:burgers}
\end{equation}
\paragraph{Related work (brief).}

PINNs \citep{lagaris1998artificial, raissi2019,karniadakis2021physics} enforce PDE structure through automatic differentiation and collocation, and are widely used as baselines for PDE surrogates. Operator-learning approaches (e.g. DeepONets and neural operators) \citep{lu2021learning, DBLP:journals/corr/abs-2108-08481} provide complementary frameworks that learn mappings from input functions to solutions. Learning low-dimensional latent representations of parametrized PDE solution manifolds has also been explored in reduced-order modeling and autoencoder-based approaches \citep{Lea_Colberg, Fresca_2022}, as well as in latent-dynamics formulations \citep{Champion_2019}. Our contribution is orthogonal: we focus on identifiable embeddings of a solution family inside a single PINN via a multihead decomposition \citep{desai2021oneshot,pellegrin2022transfer}, and on diagnostics (PCA decompositions) that summarize multiscale complexity across ensembles.

\section{Multihead PINNs as a learned basis}
\subsection{Multihead decomposition with linear heads}
We consider a \emph{family} of problems indexed by different initial conditions (ICs).
A PINN approximates $u(x,t,\nu;\mathrm{IC})$ by a neural network and minimizes a physics loss built from mean-squared residuals of \eqref{eq:burgers} over collocation points and different viscosities. This is what is called a solution bundle \cite{flamant2020solvingdifferentialequationsusing}. The IC is imposed using a hard-enforcement technique.

\paragraph{Multihead setup.}
To learn shared structures across the family, we split the network into a \emph{body} producing $n_b$ latent functions $H_j(x,t,\nu)$ and linear \emph{heads} that map the latent space to a particular solution. The final solution generated by the model can be written as follows
\begin{equation}
u(x,t,\nu;\mathrm{IC}_i)=v_i(x) + (1-e^{-t})\sum_{j=1}^{n_b} w_{ij}\,H_j(x,t,\nu),
\label{eq:linearhead}
\end{equation}
where $v_i$ is the imposed initial condition and $\{w_{ij}\}$ are the head weights for $\mathrm{IC}_i$. This exposes an explicit ``learned basis'' viewpoint: $\{H_j\}$ spans an embedding space, and the final solution is a linear combination of its components.

\paragraph{Training objective.}
Let $\mathcal{R}[u]$ be the residual of Burger's equation \ref{eq:burgers}. For each initial condition $IC_i$, we sample collocation points $(x_m,t_m, \nu_m)$ uniformly in the spacetime domain, and logarithmically in the viscosity. We minimize
\begin{equation}
\mathcal{L}=\frac{1}{n_b}\sum_{i = 1}^{n_b}\frac{1}{M}\sum_{m=1}^{M}\frac{\bigl|\mathcal{R}[u(x_m,t_m, \nu_m;\mathrm{IC}_i)]\bigr|^2}{\Lambda(\partial_xu,a,b)}
\; +\; \lambda_{\rm ortho}\,\mathcal{L}_{\rm ortho},
\label{eq:loss_total}
\end{equation}
with $\lambda_{\rm ortho} \sim 10^{-3}$ so the DE term dominates while still selecting a consistent latent space, and $n_b$ the number of heads. The last term will be explained in detail in the next subsection. The factor $\Lambda(\partial_x u,a,b)$ is a gradient-based weighting factor used to stabilize training in regions of large solution gradients. This factor weights less the regions in which the derivative of the unknown function with respect to the spatial dimension is larger. For details about the implementation we refer the reader to \cite{FerrerSanchez2024}.

\subsection{Head orthogonalization and identifiable PCA}
Embedding space functions are typically non-identifiable: different trainings can produce rotated mixtures of the same subspace.
We mitigate this by encouraging the head matrix $W=[w_{ij}]$ to be (approximately) orthonormal via the auxiliary penalty
\begin{equation}
\mathcal{L}_{\rm ortho}=\|WW^\top-I\|_F^2+\|W^\top W-I\|_F^2,
\label{eq:ortho}
\end{equation}
included in the training objective with a small weight $\lambda_{\rm ortho}$ (\ref{eq:loss_total}).
This constraint fixes the linear mixing between latent functions, rendering the embedding space identifiable. Under near-orthogonality, the principal-component decomposition of the latent space (estimated via the covariance of $\{H_j\}$ over sampled $(x,t,\nu)$) becomes stable across runs, enabling training-robust comparisons between experiments and different random initializations. This method empirically yields eigenvalues of the covariance matrix that are independent of the initial seed and allows us to extract non-degenerate information from the principal components of the latent space. Although we have applied this method to Burgers' equation, we believe that it can be generalized to other ordinary and partial differential equations.

\paragraph{How we compute PCA.}
To summarize the learned embedding, we treat the body output $H(x,t,\nu)\in\R^{n_b}$ as a random vector induced by sampling $(x,t,\nu)$ from the training domain and by mixing over the training ICs through the head weights.
Concretely, we draw a large batch of triplets $(x,t,\nu)$ and record the standardized latent vectors $\widetilde H$. The empirical covariance $\widehat\Sigma=\widetilde H\widetilde H^\top$ defines eigenvalues $\lambda_1\ge\lambda_2\ge\cdots$ and explained-variance ratios $\lambda_j/(\sum_k\lambda_k)$.
The orthogonality regularizer in \eqref{eq:ortho} makes these spectra reproducible across random initializations, so differences between different ICs ensembles can be attributed to properties intrinsic to the PDE family rather than to an arbitrary rotation of the latent basis.

\section{Experiments and results}
We train a shared body (5 fully-connected layers, width 128; $n_b=20$) with multiple linear heads, jointly over 25 viscosities $\nu\in[10^{-2},1]$ and 20 initial conditions.
We report loss curves and explained-variance spectra for two ensembles of initial data. No baseline results are shown since our approach doesn't aim to make the model converge faster. It aims into the direction of gaining interpretability on the latent space components.

\paragraph{Implementation details.}
The body takes $(x,t,\nu)$ as input, uses smooth activations (tanh), and outputs the latent vector $H$.
We use Adam with a warm-up phase during $1000$ epochs, an initial learning rate (LR) of $10^{-3}$, and a learning rate scheduler that reduces the LR by a factor of $0.985$ every $1000$ epochs. The $x$ and $t$ points are sampled from an equally spaced distribution, and $\nu$ points are logarithmically spaced. The number of points along $x$ and $t$ is $100$, whereas for the viscosity $25$ values are used. Tests were performed on an NVIDIA H100 GPU, and took approximately $  5$ days of computational time.

\subsection{Fourier initial conditions}
\begin{figure}[t]
\centering
\includegraphics[width=0.49\linewidth]{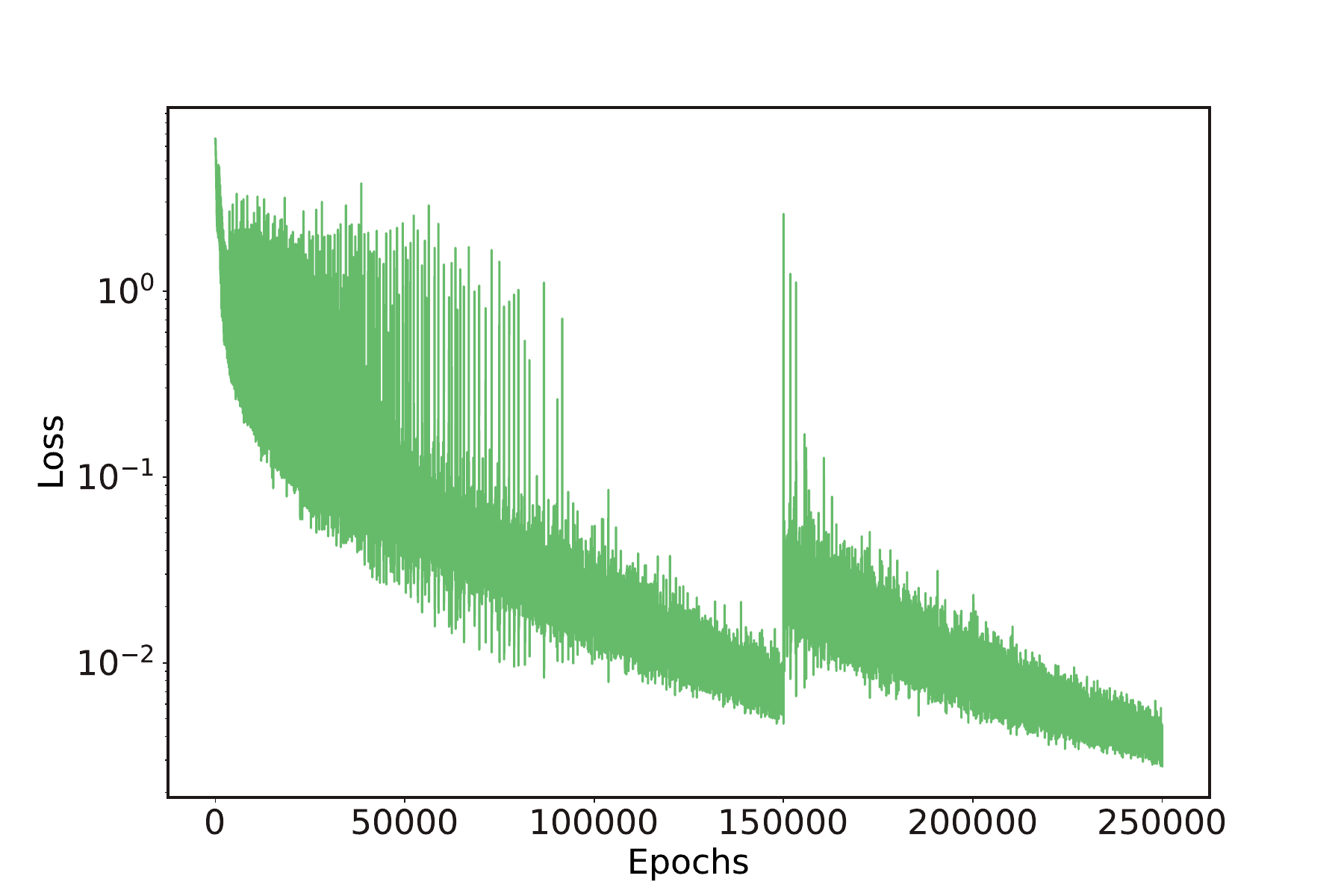}
\includegraphics[width=0.49\linewidth]{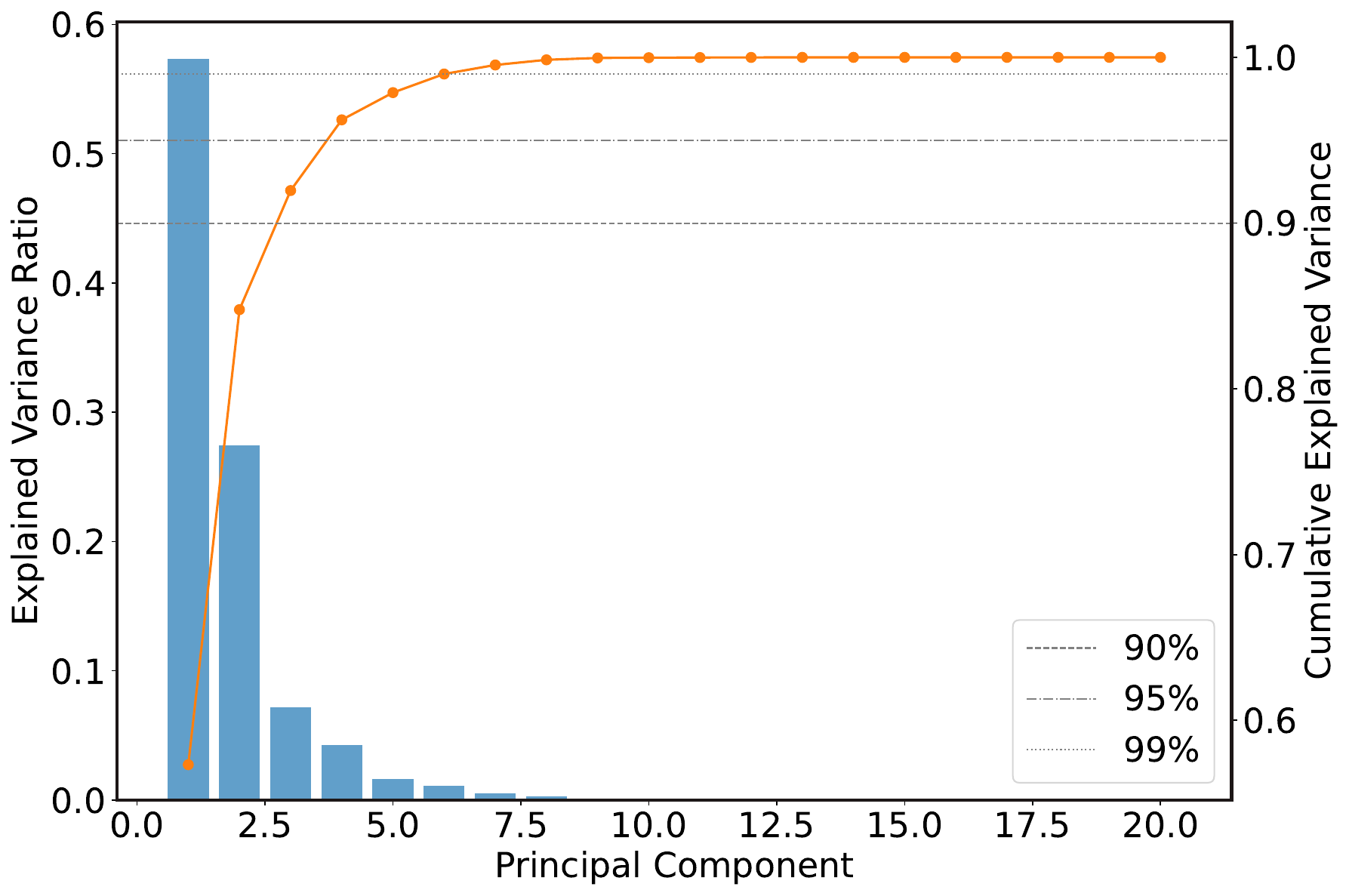}
\caption{Left: PINN PDE-residual loss versus epochs (Fourier IC ensemble). Right: Explained-variance ratio and cumulative variance of latent-space PCA (Fourier IC ensemble).}
\label{fig:fourier}
\vspace{-1ex}
\end{figure}

We construct the initial conditions as linear combinations of 10 sine and 10 cosine modes with random amplitudes and low frequencies, each modified to satisfy vanishing boundary conditions. Training converges smoothly (Figure~\ref{fig:fourier}). The latent PCA spectrum (Figure~\ref{fig:fourier}) shows rapid saturation: the first few components explain most of the variance, indicating that the solution manifold across different initial conditions admits a compact description in the learned coordinates.  The jump on the loss function is due to a change on the LR of the model, rather than something intrinsinc to the set of initial condition.

\subsection{Polynomial initial conditions}
\begin{figure}[t]
\centering
\includegraphics[width=0.49\linewidth]{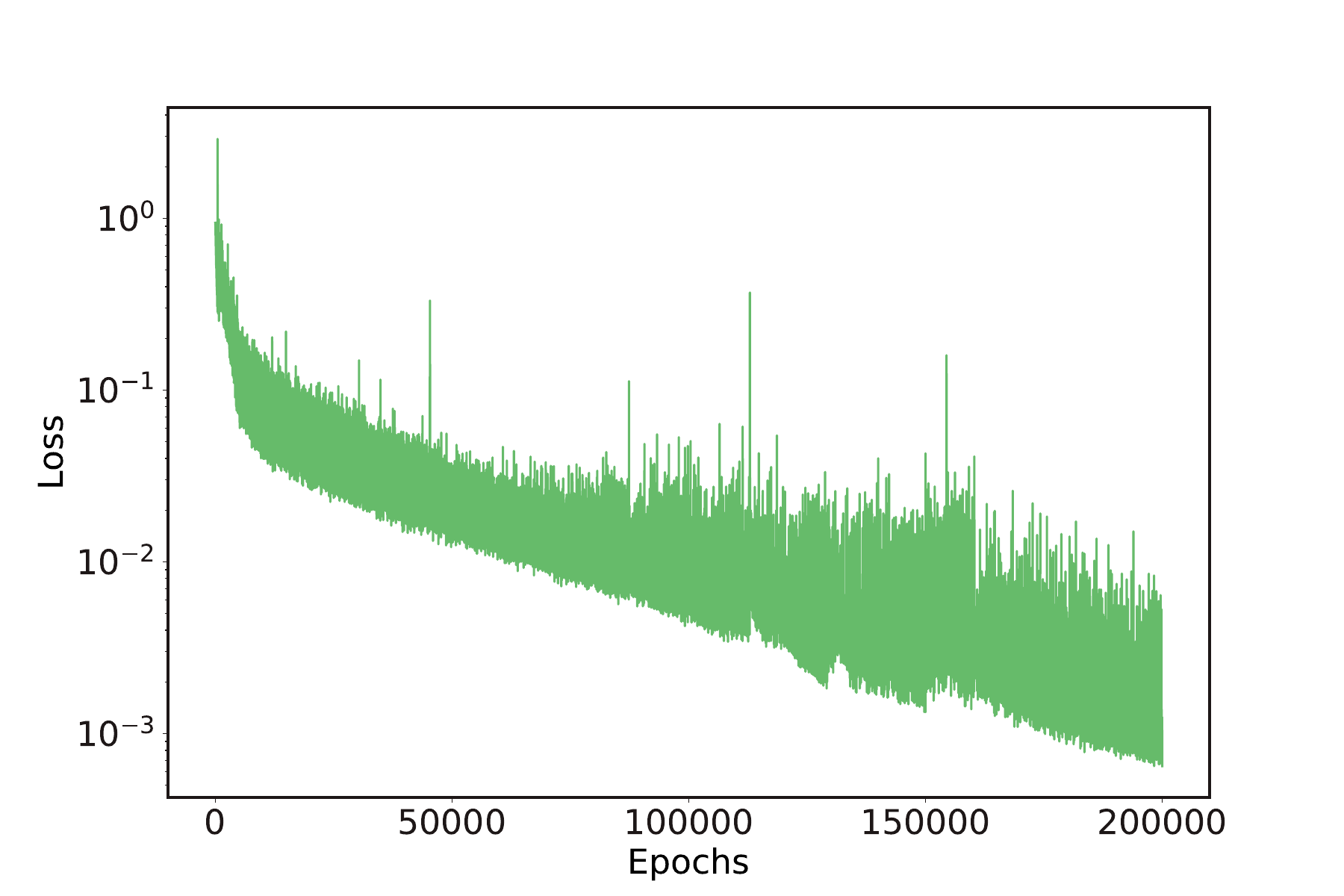}
\includegraphics[width=0.49\linewidth]{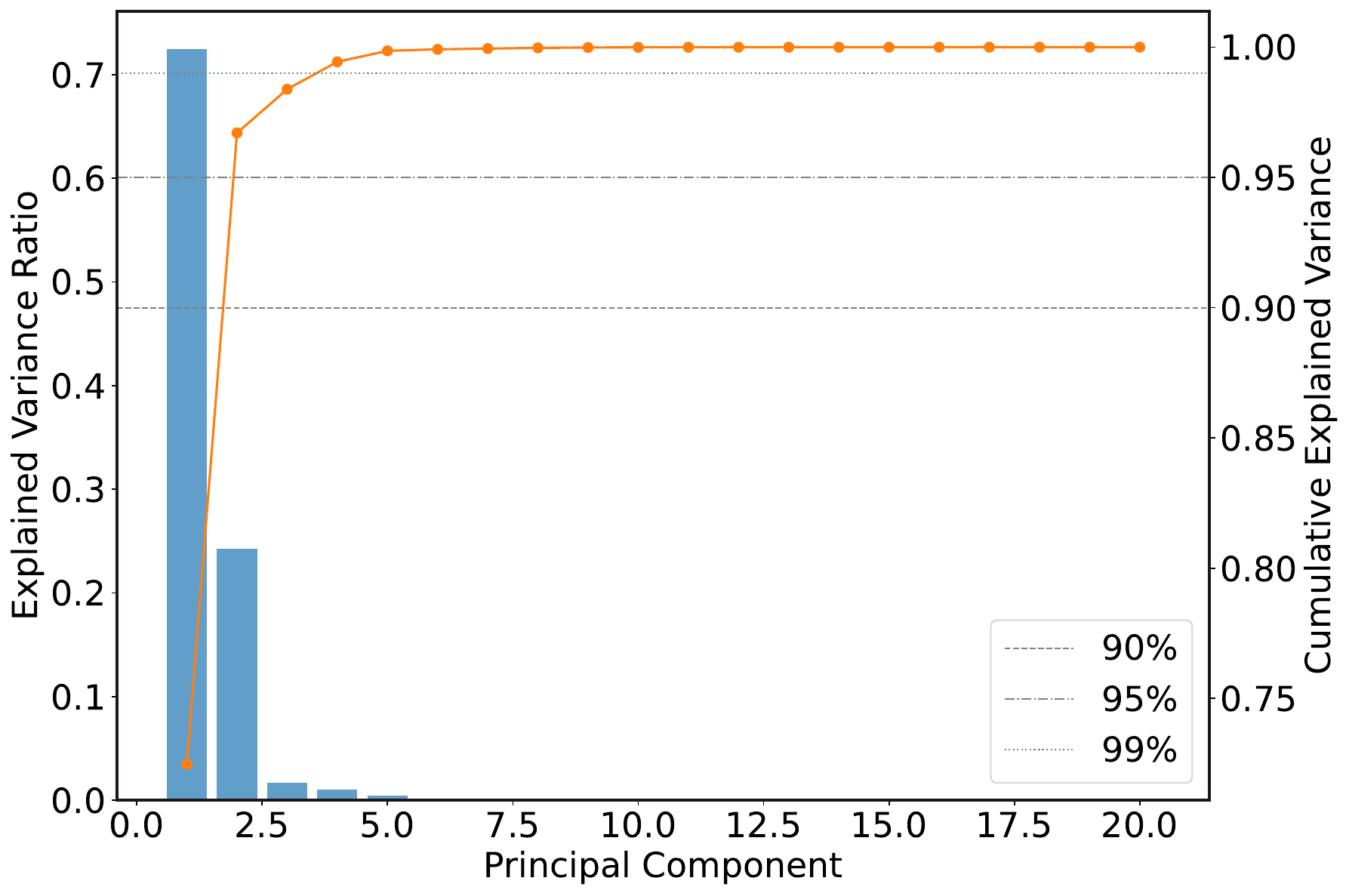}
\caption{Left: PINN PDE-residual loss versus epochs (polynomial IC ensemble). Right: Explained-variance ratio and cumulative variance of latent-space PCA (polynomial IC ensemble).}
\label{fig:poly}
\vspace{-1ex}
\end{figure}

We also generate compactly supported initial conditions as random combinations of low-degree polynomials that vanish at the edges.
Figure~\ref{fig:poly} shows the training loss.
The corresponding PCA spectrum (Figure~\ref{fig:poly}) again concentrates strongly in the leading components, and saturates even faster than in the Fourier case, consistent with a smoother, more compressible IC ensemble.

\subsection{Why embeddings matter for PDE--ML}
From the standpoint of PDE--ML, the value of this diagnostic is that it turns ``manifold complexity'' into a concrete curve: how many components are needed to capture a chosen fraction of variance.
This supports two practical workflows: (i) choosing a reduced latent dimension $r$ for downstream surrogates (e.g. learning an operator from $\mathrm{IC}$ to head weights), and (ii) stress-testing generalization by training on a larger set of collocation points and ICs, and checking whether the same leading components continue to dominate. While we focus on Burgers as a controlled proxy, the same decomposition can be applied to parameterized elliptic/parabolic PDEs and to multi-physics settings where separating coarse and fine structure is critical.

\section{Discussion: multiscale organization and an effective theory}

Across both ensembles, the embedding space organizes into a small number of hierarchical dominant components. In both cases, more than $90\%$ of the total variance can be explained by just $3$ components out of $20$. This supports an empirical \emph{effective-theory} view, in which the space of solutions spanned by different ICs can be approximated by a truncated expansion of a subset of latent space components $\{H_j\}_{j = 1,\dots, r}$. The particular form of this components is not shown in this manuscript since it is not meaningful. This is because these components are going to depend on arbitrary choices of the latent space components while training.

We find for  Burgers' equation that the most relevant components correspond to functions capturing the global structure of the solution, whereas the rest of the components tend to capture smaller-scale features. In this language, the leading PCA directions behave like an \emph{effective} basis capturing global structure\sout{s}  across the family of solutions, and additional sub-leading components act as progressively finer corrections.

Effectively, one can truncate the latent expansion at $r\ll n_b$ components to obtain a reduced model whose complexity is selected by the data-driven spectrum. 
This provides a practical model-reduction knob that is compatible with PINN training, and could improve the convergence of the model.

\section{Conclusion}
We use PINNs with a MH  setup to learn the embedding or latent space of  Burgers' equation. This  embedding space generalizes across solutions for different ICs.The heads learn how to map the latent space components to a solution with a particular IC.

We also present a method called head orthogonalization. This allows us to train the model in such a way that the latent space admits a reproducible principal component decomposition that is robust to training degeneracies.
For Burgers, a handful of components captures most of the variability across different ICs.

\paragraph{Outlook.}
For PDE--ML applications, the next step is to define a metric over the latent space manifold, in a similar way as in \cite{Tarancon-Alvarez:2025wux}. We can then construct a covariance matrix whose eigenvalues are independent of arbitrary latent space parametrization without any additional loss in expression \ref{eq:loss_total}. This will allow us to explore different future directions. The first one is identifying the reduced dimensionality of the equation with physical parameters of our model. The second consists on using a reduced model for the latent space to do transfer learning (TL) and generate different solutions on different regimes, with just a fine tuning process of the corresponding head.

It will also be interesting to apply this method to different problems as parametrized elliptic PDEs, reaction--diffusion systems, and eventually to Navier--Stokes families, where one may compare latent hierarchies with physically motivated multiscale decompositions.

\section*{LLM usage statement}

Standard AI programming tools were used during the development of the code that yielded to the results shown in this paper. These were used to increase efficiency in tasks as debugging and modifying existing code. Large language models were used to assist in editing and polishing the manuscript, but they did not contribute to the scientific content or results.

\section*{Reproducibility statement}

All results and plots can be reproduced from the equations in the article.

\bibliographystyle{iclr2026_conference.bst}

\begin{thebibliography}{13}
\providecommand{\natexlab}[1]{#1}
\providecommand{\url}[1]{\texttt{#1}}
\expandafter\ifx\csname urlstyle\endcsname\relax
  \providecommand{\doi}[1]{doi: #1}\else
  \providecommand{\doi}{doi: \begingroup \urlstyle{rm}\Url}\fi

\bibitem[Champion et~al.(2019)Champion, Lusch, Kutz, and
  Brunton]{Champion_2019}
Kathleen Champion, Bethany Lusch, J.~Nathan Kutz, and Steven~L. Brunton.
\newblock Data-driven discovery of coordinates and governing equations.
\newblock \emph{Proceedings of the National Academy of Sciences}, 116\penalty0
  (45):\penalty0 22445–22451, October 2019.
\newblock ISSN 1091-6490.
\newblock \doi{10.1073/pnas.1906995116}.
\newblock URL \url{http://dx.doi.org/10.1073/pnas.1906995116}.

\bibitem[Desai et~al.(2021)Desai, Mattheakis, Joy, Protopapas, and
  Roberts]{desai2021oneshot}
Shaan Desai, Marios Mattheakis, Hayden Joy, Pavlos Protopapas, and Stephen~J.
  Roberts.
\newblock One-shot transfer learning of physics-informed neural networks.
\newblock \emph{arXiv preprint arXiv:2110.11286}, 2021.
\newblock \url{https://arxiv.org/abs/2110.11286}.

\bibitem[Ferrer-Sanchez et~al.(2024)Ferrer-Sanchez, Martin-Guerrero,
  de~Austri-Bazan, Torres-Forne, and Font]{FerrerSanchez2024}
Antonio Ferrer-Sanchez, Jose~D. Martin-Guerrero, Roberto~Ruiz de~Austri-Bazan,
  Alejandro Torres-Forne, and Jose~A. Font.
\newblock Gradient-annihilated pinns for solving riemann problems: Application
  to relativistic hydrodynamics.
\newblock \emph{Computer Methods in Applied Mechanics and Engineering},
  424:\penalty0 116906, May 2024.
\newblock ISSN 0045-7825.
\newblock \doi{10.1016/j.cma.2024.116906}.

\bibitem[Flamant et~al.(2020)Flamant, Protopapas, and
  Sondak]{flamant2020solvingdifferentialequationsusing}
Cedric Flamant, Pavlos Protopapas, and David Sondak.
\newblock Solving differential equations using neural network solution bundles,
  2020.
\newblock URL \url{https://arxiv.org/abs/2006.14372}.

\bibitem[Fresca \& Manzoni(2022)Fresca and Manzoni]{Fresca_2022}
Stefania Fresca and Andrea Manzoni.
\newblock Pod-dl-rom: Enhancing deep learning-based reduced order models for
  nonlinear parametrized pdes by proper orthogonal decomposition.
\newblock \emph{Computer Methods in Applied Mechanics and Engineering},
  388:\penalty0 114181, January 2022.
\newblock ISSN 0045-7825.
\newblock \doi{10.1016/j.cma.2021.114181}.
\newblock URL \url{http://dx.doi.org/10.1016/j.cma.2021.114181}.

\bibitem[Karniadakis et~al.(2021)Karniadakis, Kevrekidis, Lu, Perdikaris, Wang,
  and Yang]{karniadakis2021physics}
George~Em Karniadakis, Ioannis~G Kevrekidis, Lu~Lu, Paris Perdikaris, Sifan
  Wang, and Liu Yang.
\newblock Physics-informed machine learning.
\newblock \emph{Nature Reviews Physics}, 3\penalty0 (6):\penalty0 422--440,
  2021.

\bibitem[Kovachki et~al.(2021)Kovachki, Li, Liu, Azizzadenesheli, Bhattacharya,
  Stuart, and Anandkumar]{DBLP:journals/corr/abs-2108-08481}
Nikola~B. Kovachki, Zongyi Li, Burigede Liu, Kamyar Azizzadenesheli, Kaushik
  Bhattacharya, Andrew~M. Stuart, and Anima Anandkumar.
\newblock Neural operator: Learning maps between function spaces.
\newblock \emph{CoRR}, abs/2108.08481, 2021.
\newblock URL \url{https://arxiv.org/abs/2108.08481}.

\bibitem[Lagaris et~al.(1998)Lagaris, Likas, and
  Fotiadis]{lagaris1998artificial}
Isaac~E. Lagaris, Aristidis Likas, and Dimitrios~I. Fotiadis.
\newblock Artificial neural networks for solving ordinary and partial
  differential equations.
\newblock \emph{IEEE Transactions on Neural Networks}, 9\penalty0 (5):\penalty0
  987--1000, 1998.

\bibitem[Lee \& Carlberg(2018)Lee and Carlberg]{Lea_Colberg}
Kookjin Lee and Kevin Carlberg.
\newblock Model reduction of dynamical systems on nonlinear manifolds using
  deep convolutional autoencoders.
\newblock \emph{CoRR}, abs/1812.08373, 2018.
\newblock URL \url{http://arxiv.org/abs/1812.08373}.

\bibitem[Lu et~al.(2021)Lu, Jin, Pang, Zhang, and Karniadakis]{lu2021learning}
Lu~Lu, Pengzhan Jin, Guofei Pang, Zhongqiang Zhang, and George~Em Karniadakis.
\newblock Learning nonlinear operators via deeponet based on the universal
  approximation theorem of operators.
\newblock \emph{Nature Machine Intelligence}, 3\penalty0 (3):\penalty0
  218--229, 2021.

\bibitem[Pellegrin et~al.(2022)Pellegrin, Bullwinkel, Mattheakis, and
  Protopapas]{pellegrin2022transfer}
Rapha{\"e}l Pellegrin, Blake Bullwinkel, Marios Mattheakis, and Pavlos
  Protopapas.
\newblock Transfer learning with physics-informed neural networks for efficient
  simulation of branched flows.
\newblock \emph{arXiv preprint arXiv:2211.00214}, 2022.

\bibitem[Raissi et~al.(2019)Raissi, Perdikaris, and Karniadakis]{raissi2019}
M.~Raissi, P.~Perdikaris, and G.E. Karniadakis.
\newblock Physics-informed neural networks: A deep learning framework for
  solving forward and inverse problems involving nonlinear partial differential
  equations.
\newblock \emph{Journal of Computational Physics}, 378:\penalty0 686--707,
  February 2019.
\newblock ISSN 0021-9991.
\newblock \doi{10.1016/j.jcp.2018.10.045}.

\bibitem[Taranc{\'o}n-{\'A}lvarez et~al.(2025)Taranc{\'o}n-{\'A}lvarez,
  Tejerina-P{\'e}rez, Jimenez, and Protopapas]{Tarancon-Alvarez:2025wux}
Pedro Taranc{\'o}n-{\'A}lvarez, Pablo Tejerina-P{\'e}rez, Raul Jimenez, and
  Pavlos Protopapas.
\newblock {Efficient PINNs via multi-head unimodular regularization of the
  solutions space}.
\newblock \emph{Commun. Phys.}, 8\penalty0 (1):\penalty0 335, 2025.
\newblock \doi{10.1038/s42005-025-02248-1}.

\end{thebibliography}

\end{document}